\documentstyle{article}
\newtheorem{Th}{Theorem}

\newtheorem{Lemma}{Lemma}
\newtheorem{Cor}{Corollary}
\newtheorem{Def}{Definition}
\newtheorem{Rem}{Remark}
\newcommand{\eqdef}{\stackrel{{\rm def}}{=}}

\newcommand{\dif}{\mbox{\rm d}}

\title{Dynamical and Topological methods in 
Theory of Geodesically Equivalent Metrics.}
\author{ } 
\date{}
\begin{document}
\maketitle
\vspace{0cm}

\vspace{-1.5cm}

\centerline{ Vladimir S. Matveev, 
Petar J. Topalov\footnote{Partially supported by MESC', grant number  MM-810/98.}}

\vspace{2ex}

\vspace{2ex}

Let $g=(g_{ij})$ and  $\bar g=(\bar g_{ij})$ be  $C^2$--smooth
metrics on the  same manifold $M^n$. 
\begin{Def}\hspace{-2mm}{\bf .}  The 
metrics $g$  and $\bar g$  
are {\em geodesically equivalent},
if they have  the same geodesics  (considered as unparameterized curves).
\end{Def}
This is  rather classical material. The first nontrivial examples
of geodesically equivalent metrics were constructed in 1865
by Beltrami see \cite{Bel1,Bel2}.
In  1869 Dini \cite{Dini} formulated the problem of 
local  classification of geodesically equivalent metrics, and solved it for 
 dimension two.  
 In 1896 Levi-Civita \cite{Levi-Civita}
 obtained a  local description of geodesically
equivalent metrics on manifolds of arbitrary dimension.

Many interesting results in this area were obtained by 
P. Painlev\'e, H. Weyl, E. Cartan,
P. A. Shirokov, S. Kobayashi, N. S. Sinyukov,
A. Z.   Petrov, P. Venzi, J. Mike\v{s}, A. V.  A\-mi\-no\-va, see
\cite{Amin, Mikes1}  for references.

The main tool they used was tensor analysis, and the most results
were local. 
The approach we would like to suggest is more global; in particular 
it helps to find an answer for the 
following questions: 
\begin{itemize}
\item  What closed manifolds admit geodesically equivalent metrics? \cite{Mikes}
\item How many metrics are there geodesically equivalent to a given one?  \cite{Mikes1}
\end{itemize}
Our approach is based on the following theorem.

Denote by $G:TM^n\to TM^n$ the fiberwise-linear mapping given by the tensor
$g^{-1}\bar g=(g^{i\alpha}\bar g_{\alpha j})$. 
In invariant terms, for  any   $x_0\in M^n$  the
restriction of the mapping
$G$ to the tangent space  $T_{x_0}M^n$ is the
 linear transformation
of  $T_{x_0}M^n$
 such that for any  vectors
$\xi, \ \nu\in T_{x_0}M^n$
the
dot product $g(G(\xi), \nu)$ of the vectors $G(\xi)$ and $\nu$
in the metric $g$
  is equal to
the
dot product $\bar g(\xi, \nu)$ of the vectors $\xi$ and $\nu$
in the metric $\bar g$.
Consider the characteristic polynomial  
$det(G- \mu E)=c_0\mu^n+c_1\mu^{n-1} + ... + c_n$. The
coefficients $c_1, .. , c_n$ are smooth functions on the 
manifold $M^n$, and $c_0\equiv (-1)^n$. 
Consider the fiberwise-linear mappings
$$ S_0, S_1,  ... ,S_{n-1}:TM^n\to TM^n $$ 
given by the 
 general formula  
$$
S_k\eqdef\left(\frac{det(g)}{det(\bar g)}\right)^{\frac{k+2}{n+1}}\sum_{i=0}^kc_i G^{k-i+1}
.$$

Consider the functions $I_k:T^*M^n\to R$, $k=0,...,n-1$,
 given 
by the formulae
\begin{equation}\label{12345}
 I_k(x,p)\eqdef  g^{\alpha i}(S_k)_\alpha^j p_ip_j.
\end{equation}
In invariant terms the functions $I_k$ are as follows.
Let us   identify canonically
the tangent and the cotangent bundles of
 $M^n$ by the
metric $g$.
Then the value of the function
$I_k$ on a vector $\xi\in T_{x}M^n\cong T_x^*M^n$ is given by
\begin{equation}\label{I}
I_k(x, \xi)\eqdef
g(S_k(\xi),\xi).\
\end{equation}
 Consider the standard symplectic form on $T^*M^n$.
 By geodesic flow of the metric $g$  we mean  the Hamiltonian system
 on $T^*M^n$ with the Hamiltonian
 $H(\xi)=\frac{1}{2}g(\xi,\xi)$.

\begin{Rem}\hspace{-2mm}{\bf .}
The function  ${I}_{n-1}$ is equal to minus twice Hamiltonian of the
geodesic flow of the metric $g$.
 \end{Rem}

\begin{Th} \hspace{-2mm}{\bf .}\label{classical}
If the metrics  $g$ and $\bar g$ on $M^n$  
are  geodesically equivalent 
then   
the functions 
$I_k$
 pairwise  commute.
 In particular, they are
 integrals in involution
 of the geodesic flow of the metric $g$.
\end{Th} 
This theorem allows us to apply the theory of integrable geodesic flows to the
 theory of geodesically equivalent 
metrics and vice versa. 
First of all, the following theorem explaines
  when the integrals are functionally independent.

\begin{Th}[\cite{TM}]\hspace{-2mm}{\bf .}
\label{123}
   Let metrics $g, \ \bar g$ be geodesically equivalent.   
If   at a point  $x_0\in M^n$ 
the number of different 
roots  of the polynomial  $det(G- \mu E)$ is equal to
 $n_1$,  then at almost all points of $T_{x_0}M^n$ 
the dimension of the linear space generated by $dI_0, dI_1, ... ,dI_{n-1}$ is greater  or equal than $n_1$.   

In particular, 
if the characteristic polynomial $det(G- \mu E)$ has no
multiple roots at a point $x_0\in M^n$
then the 
functions $I_0,I_1,..., I_{n-1}$ are
functionally independent almost everywhere in $TU(x_0)$, where 
$U(x_0)$ is a sufficiently  small neighborhood of $x_0$.

If   at every point of  a
neighborhood $U(x_0)$ of a  point $x_0\in M^n$ 
the number of different 
eigenvalues of the polynomial  $det(G- \mu E)$ is less or equal than 
$n_1$, 
then at all points of $TU(x_0)$  
the dimension of the linear space generated by $dI_0, dI_1, ... ,dI_{n-1}$ is less or equal than $n_1$.   
\end{Th}
The metrics $g, \bar g$ are    {\em strictly  non-proportional} at a point $x_0\in M^n$ when
 the characteristic polynomial $det(G- \mu E)_{|x_0}$
 has no multiple roots.

\begin{Cor}\hspace{-2mm}{\bf .} \label{connected}
Suppose $M^n$ is  connected. 
Let metrics $g, \ \bar g$ on 
$M^n$ 
be geodesically equivalent and 
 strictly  non-proportional at least 
at one point on  $M^n$. 
Then the metrics  are strictly non-proportional almost everywhere. 
 \end{Cor}
More generally, suppose    $M^n$ is    connected and
 the  metrics $g, \ \bar g$ on 
$M^n$ 
are  geodesically equivalent.  If 
 at every point of  a
neighborhood $U\subset  M^n$ 
the number of different 
eigenvalues of the polynomial  $det(G- \mu E)$ is less or equal than $n_1$, 
 then at every point of  
$M^n$ the number of different 
eigenvalues of the polynomial  $det(G- \mu E)$ is less or equal than $n_1$.

Corollary~\ref{connected}  follows from the following observation.
If we have an integrable Hamiltonian systems then the  dimension of the linear space generated by 
differentials of the integrals is constant along each orbit. 

Proof of Corollary~\ref{connected}.   
Identify canonically
the tangent and the cotangent bundles of
 $M^n$ by the
metric $g$.
Take a geodesic $\gamma:R\to M^n$ 
and consider the points $x_0=\gamma(0), \ x_1=\gamma(1)\in M^n. $
Suppose that the metric $\bar g$ is geodesically equivalent to the metric
  $\bar g$ and is  strictly non-proportional with $g$ at the point $x_0$. Let 
us prove that  in  each neighborhood of the point $x_{1}$ there
exists  a point $x$ such that  the metrics are strictly non-proportional  
at $x$.

The geodesic orbit of the geodesic $\gamma$ is the curve 
 $(\gamma,\dot\gamma): R\to TM^n\cong T^*M^n$
 (assuming $\dot\gamma=\dot\gamma(t)\in T_{\gamma(t)}M^n$
 is the velocity vector of $\gamma$.)
Take  a small neighborhood 
 $O^{2n}\subset  TM^n$  
of the point $(x_0, \dot\gamma(0))\in  TM^n$. 
Consider the union of all 
points of  
all geodesic orbits that start  from the points of the neighborhood. 
Denote it by $U$.  Since the  metrics are  strictly non-proportional at $x_0$, 
the differentials $dI_0,dI_1, ... ,dI_{n-1}$ 
of the integrals $I_1, I_2, ...,I_{n-1}$ are linearly  
independent at almost all points 
of  $U$. Evidently the set $U$ contains a neighborhood of the 
point $(x_1, \dot\gamma(1))$. The integrals 
 $I_0,I_1,..., I_{n-1}$ are then  functionally independent at almost all 
points of a  neighborhood of the 
point $(x_1, \dot\gamma(1))$. Using Theorem~{\ref{123}},
 we have 
that
in each neighborhood of  $x_1$ there exists a point where
the metrics are 
strictly non-proportional.  Since each two points of a connected
manifold can be joint by a sequence of geodesical segments,
in each neighborhood of  an arbitrary point of $M^n$
 there exist points where
the metrics are 
strictly non-proportional, q. e. d.

Let us combine Theorems~\ref{classical},\ref{123}
with known results
in the theory of integrable geodesic flows.

\begin{Th}[Taimanov, \cite{Taiman}]\hspace{-2mm}{\bf .} \label{Taim} 
If  a real-analytic closed manifold $M^n$ with a real-analytic
metric satisfies at least one of the conditions:
\begin{itemize}
\item[a)] $\pi_1(M^n)$ is not almost commutative
\item[b)] $dim H_1(M^n;{\bf Q})>dim M^n$, 
\end{itemize}
then the geodesic flow on $M^n$ is not analytically integrable. 
\end{Th}

\begin{Cor}[\cite{TM}]\hspace{-2mm}{\bf .} \label{T} Let $M^n$ be a closed real-analytic 
manifold supplied with
  two real-analytic metrics $g,\bar g$ such that the metrics $g,\bar g$ are 
 geodesically equivalent and 
                 strictly non-proportional at  least at one point.     
Then  the fundamental group 
$\pi_1(M^n)$ of the manifold $M^n$ contains  a commutative   subgroup of finite 
index and the dimension of the homology  group $H_1(M^n; {\bf Q})$  is 
 no greater than $n$.  
\end{Cor}
Easy to see that the integrals $I_k$ are quadratic in momenta.

\begin{Th}[\cite{Kol}]\hspace{-2mm}{\bf .} Let $g$ be a
Riemannian metric on a closed
surface $M^2$. Suppose that  the geodesic flow of $g$ admits an integral
that is
quadratic in momenta and functionally independent of the Hamiltonian. 
Then  the surface $M^2$ is homeomorphic either  to the torus or to
 the sphere or to the Klein bottle or to the projective plane.
\end{Th}

\begin{Cor}[\cite{TM}]\hspace{-2mm}{\bf .}{\label{Kol}} Let metrics $g, \bar g$ on a closed 
surface of 
negative Euler characteristic be 
geodesically equivalent.
 Then $g=C  \bar g$, where $C$ is a constant.  
\end{Cor}
We say that the Hamiltonian and an  integral of  a geodesic flow
are {\it proportional} at
a point  $x\in M^n$,  if  for a constant $c$ for each
point
$\xi\in T^*_xM^n$ we have $I(\xi)=cH(\xi)$.

Recall that a vector field on $M^n$  is 
{\it Killing} (with respect to a metric), 
 if the flow of the field 
 preserves the metric. 

\begin{Th}[\cite{Kol}]\hspace{-2mm}{\bf .}\label{sphere}
Let the geodesic flow of a  metric $g$
 on the sphere $S^2$ admits an integral $I$ quadratic in velocities 
and functionally independent of the Hamiltonian $H$ of the geodesic flow.
Then there are 
only three possibilities. 
\begin{itemize}
\item[1.] The Hamiltonian and the integral
are proportional  at   exactly   two points.
\item[2.] The Hamiltonian and the integral are proportional
at  exactly   four points. 
\item[3.] The    Hamiltonian and the integral   are completely 
proportional, i.e. $c H= I$,
 where $c$ is a constant.
\end{itemize}
In the first case the  metric admits a Killing vector field. 
\end{Th}
Easy to see that the metrics $g, \ \bar g$
 are proportional at a point if and only if the Hamiltonian of the 
 geodesic flow of $g$ and the integral $I_0$ are proportional at 
the point.

\begin{Cor}[\cite{TM}]\hspace{-2mm}{\bf .}\label{sphere1}
Let metrics $g$, 
$\bar g$ on the sphere $S^2$  be geodesically equivalent. 
Then there are  only three possibilities. 
\begin{itemize}
\item[1.] The metrics are proportional  at   exactly   two points.
\item[2.] The metrics are proportional  at  exactly   four points. 
\item[3.] The metrics are completely proportional, i.e. $g=C  \bar g$,
 where $C$ is a positive constant.   
\end{itemize}
In the first case the  metrics admit a Killing vector field. 
\end{Cor}

\begin{Cor}\hspace{-2mm}{\bf .}\label{killing}
Let   $g,\bar g$ be geodesically equivalent metrics
   on $M^n$.
   If
   $g$ admits  a non-trivial  Killing vector field
 then $\bar g$ also admits
 a non-trivial Killing vector
field. 
\end{Cor}
\noindent { Proof.}   Let $g$, $\bar g$ be geodesically
equivalent metrics on $M^n$.
By  Noether's  theorem,
if a metric admits a (non-trivial) Killing
vector field, 
then the geodesic flow of the metric admits a (non-trivial) 
 integral, linear in velocities,  and vice versa. 
Then it is sufficient to prove that, given
an integral linear in velocities
for the geodesic flow of $\bar g$,  we can  construct an integral linear
in velocities for the geodesic flow of
$ g$.
Since the metrics are geodesically   equivalent, the fucntion
$\phi:T^*M^n\cong TM^n\to TM^n\cong T^*M^n$
given by
\begin{equation}\label{orbit}
\phi(x,\xi)=\left(x,\frac{\|\xi\|_g}{\|\xi\|_{\bar g}}\xi\right)
\end{equation}
takes the orbits of the geodesic flow of $g$ to the orbits of
the geodesic flow of $\bar  g$.
    Suppose the
 function 
$$F_1=\sum^n_{i=1}a_i(x)\xi^i$$
is constant on 
the trajectories of the geodesic flow of the metric
$\bar g$. Then  the function
$$\phi^*F_1=
\frac{||\xi||_{g}}{||\xi||_{\bar g}}\sum^n_{i=1}a_i(x)\xi^i$$ is 
constant on the trajectories of the geodesic flow of the metric $g$. Since the function 
$
I_0=\left(\frac{{\rm det}(g)}{{\rm det}(\bar g)}\right)^{\frac{2}{n+1}}\bar g(\xi, \xi)
$
is an integral  of the geodesic flow of the metric $g$, and since the function 
$||\xi||_{g}=\sqrt{g(\xi,\xi)}$ is also an integral 
of the geodesic flow of the metric $g$, then the 
function  
$$\frac{\sqrt{g(\xi,\xi)}}{\sqrt{I_0}}\phi^*F_1=
\left(\frac{{\rm det}(g)}{{\rm det}(\bar g)}\right)^{\frac{1}{n+1}}\sum^n_{i=1}a_i(x)\xi^i, 
$$ linear in velocities,  
is also an integral 
of the geodesic flow of the metric $g$, q.~e.~d. 

         \vspace{3mm}

Does the integrability of a geodesic flow implies the existence 
of geodesically equivalent metric?
Theorem~\ref{classical} shows that 
 we should   have   integrals
in the form (\ref{12345}).

\begin{Th}\hspace{-2mm}{\bf .}\label{inverse}
Let $g, \ \bar g$  be metrics on $M^n$. Consider the functions 
$I_0,I_1,  ... ,I_{n-1}$. Let they be functionally 
 independent almost everywhere and let they  commute. 
 Then the metrics $g,\  \bar g$ are geodesically equivalent.   
\end{Th}

\vspace{2ex}

\noindent {\bf Dynamical Background:}

\noindent We suggest a construction that, given an 
 orbital  diffeomorphism between   
 two Hamiltonian systems,  produces integrals of them.  

  Let  $v$ and $\bar v$ be Hamiltonian systems
 on  
 symplectic manifolds $(N,\omega)$ and $(\bar N,\bar\omega)$
 with Hamiltonians $H$ and
$\bar H$,  respectively.

Consider the 
isoenergy surfaces  
$$
Q\eqdef\left\{x\in N: H(x)=h\right\}, \ \ 
\bar Q\eqdef\left\{x\in \bar  N: \bar H(x)=\bar h\right\},
$$
\noindent where $h$ and $\bar h$ are regular values of the functions 
$H$, $\bar H$,  respectively.

  \begin{Def}\hspace{-2mm}{\bf .}\label{TrAutDef}
A diffeomorphism $\phi : Q\longrightarrow \bar Q$ 
 is said to be  {\em
orbital}, 
 if it 
 takes the orbits  of the system $v$ to
the orbits of the system $\bar v$.
\end{Def}
 Given orbital diffeomorphism, we can invariantly construct integrals.
 Denote by $\sigma$, $ \bar \sigma$
the restrictions of the forms $\omega, \bar \omega$  to $Q$, $\bar Q$  
respectively.    
Consider the form $\phi^*\bar \sigma$ on $Q$.

\begin{Lemma}[Topalov, \cite{TensInv}]\hspace{-2mm}{\bf .}\label{TensInv}
The flow  $v$ preserves the
 form  $\phi^{*}\bar\sigma$.
\end{Lemma}
\noindent { Proof.}
The  Lie derivative $L_v$ of the form $\phi^*\bar \sigma$ along the vector field $v$ satisfies  
$$
L_v\phi^{*}\bar\sigma=\dif\left[{\imath}_v\phi^{*}\bar\sigma\right]+{\imath}_v
\dif\left[\phi^{*}\bar\sigma\right]. 
$$
On the right-hand  side  both terms vanish.
\noindent Since the form $\bar\omega$ is closed, the form $\bar \sigma$ is also closed and
$\dif\left[\phi^{*}\bar\sigma\right]=\phi^{*}(\dif\bar\sigma)=0$. 
Since the diffeomorphism takes the orbits to orbits, it takes the kernel of the form $\sigma$  to the kernel of the form $\bar \sigma$, so that 
${\imath}_v\phi^{*}\bar\sigma$ is equal to zero, q. e. d.

It is obvious  that the kernels of the forms $\sigma$ and
$\phi^{*}\bar\sigma$ coincide (in the space ${\cal T}_xQ$ at each point $x\in Q$) 
with the linear  span of
the  vector  $v$. Therefore these forms induce two
non-degenerate
tensor fields on the quotient bundle 
${\cal T}Q/\langle v\rangle$. We shall
denote the corresponding forms on ${\cal T}Q/\langle v\rangle$ also by the  letters $\sigma, \bar\sigma$.

\begin{Lemma}\hspace{-2mm}{\bf .}\label{Din}
The characteristic polynomial of
$(\sigma)^{-1}(\phi^{*}\bar\sigma)$ on ${\cal T}Q/\langle v\rangle$
is preserved by the flow $v$.
\end{Lemma}

\noindent{Proof.} Since
 the flow $v$ preserves the 
Hamiltonian  $H$ and the form $\omega$,   
the flow $v$ preserves 
the form $\sigma$. 
Since the  flow  $v$ preserves 
both forms, it preserves  
the  characteristic
polynomial of
$(\sigma)^{-1}({\phi^{*}\bar\sigma})$, q.~e.~d.

Since both forms are skew-symmetric, each root of the  characteristic
 polynomial 
$(\sigma)^{-1}(\phi^{*}\bar\sigma)$ has an even multiplicity. Then  the characteristic
 polynomial   is the  square of a polynomial 
$\delta^{n-1}(t)$ of degree $n-1$. Hence  
the polynomial $\delta^{n-1}(t)$ is also preserved by the flow $v$.
Therefore the coefficients of the polynomial $\delta^{n-1}(t)$ are integrals of the system $v$.

Geodesic flows of geodesically equivalent metrics can be considered 
as orbitally equivalent systems. The manifold 
$N=\bar N=TM^n$, the forms $\omega$, $\bar \omega$ are given by 
$$
\omega \eqdef\dif[g_{ij}\xi^jdx^i], \  \ \bar\omega \eqdef\dif[\bar g_{ij}\xi^jdx^i]
$$
and the orbital diffeomorphism  
$\phi$ is given by (\ref{orbit}).
Direct calculations give us the formulae for the integrals 
$I_k$ from Theorem~1.

\vspace{3ex}

\noindent{\bf Are there interesting examples of geodesically equivalent metrics on
closed manifolds?}

\noindent 
The following theorem (essentially due to N. S. Sinjukov \cite{Sin}) gives us 
 a construction that, given a pair of geodesically equivalent metrics, produces
another  pair of geodesically equivalent metrics. Starting from the metric 
of constant curvature on the sphere, we obtain the metric 
of the  ellipsoid and the metric of the Poisson sphere. 

Let $g, \bar g$ be Riemannian metrics on  $M^n$.
Consider the fiberwise-linear mapping   $B:TM^n\to TM^n$ given by
$B_j^i=\left(\frac{det(\bar g)}{det (g)}\right)^{\frac{1}{n+1}}{\bar g^{i\alpha}}g_{\alpha j}$.
 By definition,  let us  put metric $g_B$ equal to  $g_{ip}B^p_j$ and put
metric $\bar g_B$ equal to  $\bar g_{ip}B^p_j$.
In invariant terms, the dot product  $g_B(\xi, \nu)$
 of  arbitrary
 vectors $\xi, \nu \in T_{x_0}M^n$ is equal to $g(B\xi, \nu)$ and
the dot product  $\bar g_B(\xi, \nu)$
  is equal to $\bar g(B\xi, \nu)$.
 Evidently,  the restriction of  $B$  to any tangent space $T_{x_0}M^n$
 is self-adjoint with respect to the
 metrics $g$ and $\bar g$ and the metrics  $g_B,\  \bar g_B$ are well-definite.

\begin{Th}\hspace{-2mm}{\bf .}\label{Sinjukov}
The metrics  $g$   and  $\bar g$ are geodesically equivalent if and only if 
the metrics 
 $g_B$ and   $\bar g_B$ are geodesically equivalent.
\end{Th}
Evidently,  the metrics $g$  and  $ \bar g$ are
strictly  non-proportional at a point $x\in M^n$ 
 if and only if
the metrics  $g_B$ and $  \bar g_B$ are strictly non-proportional
at the  point $x$.

Thus if we have a pair of geodesically 
equivalent metrics $g$, $\bar g$, 
then we can construct the other 
pair of geodesically equivalent metrics  $g_B, \ \bar g_B$.
We can apply the construction
 once more, the result is another pair of geodesically equivalent metrics. 
It is natural to denote it by $g_{B^2},$    $\bar g_{B^2}$ since
this pair is given by 
$$
g_{B^2}(\xi,  \nu)=g(B^2 \xi, \nu), \ \   \bar g_{B^2}(\xi,  \nu)=\bar g(B^2 \xi, \nu).$$
We can go in other direction and   consider the metrics 
$
g_{B^{-1}},\  \bar g_{B^{-1}}$ given by
$$g_{B^{-1}}(\xi,  \nu)=g(B^{-1} \xi, \nu), \ \   \bar g_{B^{-1}}(\xi,  \nu)=\bar g(B^{-1} \xi, \nu).$$
 They  are geodesically equivalent as well.

To start the process, we need a
 pair of geodesically equivalent metrics $g$, $\bar g$. 
We take the following one (obtained by E. Beltrami \cite{Bel1}, \cite{Bel2}).
 The metric $g$ is the restriction of the Euclidean metric
 $$(dx^1)^2+(dx^2)^2+...+(dx^{n+1})^2$$ to the standard sphere 
$$
S^n=\{(x^1, x^2, ... ,x^{n+1})\in R^{n+1} : (x^1)^2+(x^2)^2+...+(x^{n+1})^2=1\}. 
$$

The metric $\bar g$ is the pull-back $l^* g$,
where the diffeomorphism
$l:S^n\to S^n$ is given by 
$$
l(x)\eqdef \frac{Ax}{\| Ax \|},
  $$ 
where $A$ is an arbitrary non-degenerate linear 
transformation of  $R^{n+1}$, $Ax$ means 'the transformation 
$A$ applied to the vector 
$x=(x^1, x^2, ..., x^{n+1})$,'
and  
$\| x\|$ is the standard norm $\sqrt{(x^1)^2+(x^2)^2+...+(x^{n+1})^2}$. 

Easy to see that the mapping $l$
preserves the geodesics of the sphere. More precisely,
the geodesics of the sphere are intersections of the planes,  which
  go through 
the origin, with  the sphere. The linear transformation
 $A$ takes
 planes to planes, therefore the mapping $l$ takes 
geodesics to geodesics. Then the standard 
 metric $g$ of the sphere and the metric
$l^*g$ are geodesically equivalent.

For these  
geodesically equivalent
metrics 
 the metric $g_{B}$ is
 the metric of an
 ellipsoid, and the metric $\bar g_{B^2}$ is the metric of a
Poisson sphere. 
By varying the linear operator $A$, we can obtain metrics
of all  possible ellipsoids  and all possible  Poisson spheres.

Recall that the metric of the
ellipsoid 
$$E\eqdef \left\{(x^1,x^2, ... ,x^{n+1})\in  R^{n+1} :
\frac{(x^1)^2}{a_1}+\frac{(x^2)^2}{a_2}+...+\frac{(x^{n+1})^2}{a_{n+1}}=1\right\}$$
is the restriction of the metric 
$$(dx^1)^2+(dx^2)^2+...+(dx^{n+1})^2$$ 
to the ellipsoid $E$.  
By the metric of the Poisson sphere we, following  \cite{Brailov}, mean  
the restriction of the metric
$$
\frac{1}{\frac{(x^1)^2}{a_1^2}+\frac{(x^2)^2}{a_2^2}+...+\frac{(x^{n+1})^2}{a_{n+1}^2}}\left((dx^1)^2+(dx^2)^2+...+(dx^{n+1})^2\right)
$$
to the ellipsoid $E$.

The metric of Poisson sphere has the following mechanical sense. 
Consider the free motion of (n+1)-dimension rigid body in (n+1)-dimensional  
space around a fixed point. 
The configuration space of the corresponding Hamiltonian system is
$SO(n+1)$,  and the corresponding Hamiltonian is left-invariant.
 Consider the embedding $SO(n)$ into $SO(n+1)$ as
the stabilizer of a
vector $v\in R^{n+1}$. Consider the action of the group $SO(n)$ on
$SO(n+1)$ by left translations. The Hamiltonian of 
the motion is evidently invariant modulo this action, 
and the reduced system 
on $TSO(n+1)/{SO(n)}\cong TS^n$ is the geodesic flow of the (appropriate)
Poisson 
metric.

\begin{Th}[\cite{TM},
independently obtained by S. Tabachnikov
\cite{Tab}]\hspace{-2mm}{\bf .}\label{ellipsoid}
 The restriction of the Euclidean metric 
 $$(dx^1)^2+(dx^2)^2+...+(dx^{n+1})^2$$ to the standard ellipsoid
$$
E^n=\left\{(x^1, x^2, ... ,x^{n+1})\in R^{n+1} : \frac{(x^1)^2}{a_1}+\frac{(x^2)^2}{a_2}+...+\frac{(x^{n+1})^2}{a_{n+1}}=1
\right\}
$$  
is geodesically equivalent to the restriction of the metric 
\begin{equation}\label{8}
\frac{1}{
\left(\frac{x^1}{a_1}\right)^2+\left(\frac{x^2}{a_2}\right)^2+...
+\left(\frac{x^{n+1}}{a_{n+1}}\right)^2
        }
\left(\frac{(dx^1)^2}{a_1}+\frac{(dx^2)^2}{a_2}+...+\frac{(dx^{n+1})^2}{a_{n+1}}\right)
\end{equation}
to the same ellipsoid. 
\end{Th}

\begin{Th}[Topalov, \cite{Top}]\hspace{-2mm}{\bf .} \label{poisson}
 The restriction of the metric 
\begin{equation}\label{9}
\frac{1}{\frac{(x^1)^2}{a_1^2}+\frac{(x^2)^2}{a_2^2}+...+\frac{(x^{n+1})^2}{a_{n+1}^2}}\left((dx^1)^2+(dx^2)^2+...+(dx^{n+1})^2\right)
\end{equation}
to the ellipsoid 
  $$
E^n=\left\{(x^1, x^2, ... ,x^{n+1})\in R^{n+1} : \frac{(x^1)^2}{a_1}+\frac{(x^2)^2}{a_2}+...+\frac{(x^{n+1})^2}{a_{n+1}}=1\right\}
$$ 
is geodesically equivalent to 
the restriction of the metric 
\begin{equation}\label{10}
  a_1(dx^1)^2+a_2(dx^2)^2+ ...
   +a_{n+1}(dx^{n+1})^2-(x^1dx^1+x^2dx^2 +...+x^{n+1}dx^{n+1})^2
\end{equation}
to the same ellipsoid. 
\end{Th}

\vspace{2ex} 

\noindent{\bf Quantum integrability of the Beltrami-Laplace operator for geodesically equivalent metrics:}

\noindent Let $g$, $\bar g$ be Riemannian metrics on $M^n$.
Consider the Beltrami-Laplace operator
 $$
 \Delta:=
 -\frac{1}{\sqrt{det(g)}}\frac{\partial }{\partial x^i}
\sqrt{det(g)} g^{ij}
 \frac{\partial }{\partial x^j}, $$
 where $det(g)$ denotes the  determinant of the matrix
  corresponding to the
  metric $g$.

Consider the operators 
\begin{equation}\label{I-quantum}
{\cal I}_0,  \ {\cal I}_1, ... ,{\cal I}_{n-1}:C^2(M^n)\to C^0(M^n)
\end{equation}
given 
 by the
 general formula
\begin{equation} \label{I-quantum1}
{\cal I}_k:=
 \frac{1}{\sqrt{det(g)}}\frac{\partial }{\partial x^i} (S_k)^i_\alpha
\sqrt{det(g)} g^{\alpha j}
 \frac{\partial }{\partial x^j}.
\end{equation}

\begin{Rem}\hspace{-2mm}{\bf .}
The operator ${\cal I}_{n-1}$ is exactly the operator $\Delta$.
 \end{Rem}

\begin{Th}\hspace{-2mm}{\bf .}\label{main}
If the metrics  $g$ and $\bar g$ on $M^n$  
are  geodesically equivalent  then 
the operators  
${\cal I}_k$ pairwise  commute. In particular they 
 commute with the Beltrami-Laplace operator 
$\Delta$.  
If the manifold $M^n$ is closed then the operators ${\cal I}_k$ are 
self-adjoint.  
\end{Th}

\begin{Cor} \hspace{-2mm}{\bf .}
   Suppose $M^n$ is connected. 
  Metrics $g, \ \bar g$ on $M^n$ 
  are geodesically equivalent and strictly non-proportional 
  at least at  one point of  $M^n$ if and  only if the     
 operators  ${\cal I}_0, ... ,{\cal  I}_{n-1}$ are 
 linearly  independent.   
\end{Cor}

\begin{Cor}\hspace{-2mm}{\bf .}
Suppose the manifold $M^n$ is connected. Let  $g, \ \bar g$ be
 metrics on it.  
Let the operators ${\cal I}_k$ commute and let they be 
linearly independent.
Then
the metrics $g, \bar g$ are geodesically equivalent. 
\end{Cor}
So that if the manifold is closed and if the
metrics $g, \ \bar g$ on it
  are geodesically equivalent and strictly non-proportional 
  at least at  one point then we have the complete quantum integrability of
 the Beltrami-Laplace operator of the metric $g$.

Quantum integrability means that 
there exists a countable   basis 
 $$\Phi=\{\phi_1, \phi_2, ..., \phi_m, ...\}$$
 of the space $L_2(M^n)$
 such that each $\phi_m$ is an eigenfunction of each  operator  
${\cal I}_k$.

 Moreover, in our case the variables can be separated. More precisely, 
take any function $\phi$ from the basis $\Phi$.
 Since 
$\phi$ is an eigenfunction of each operator ${\cal I}_k$,
we have that $\phi$ is a solution of the system of $n$ partial 
differential equations 
\begin{equation}\label{bl}
{\cal I}_k\phi=\lambda_k\phi, \  \ k=0,1,...,n-1.
\end{equation}
The separation of variables means that 
 in a neighborhood    
of almost any point there exist coordinates
 $(x^0,x^1, ... ,x^{n-1})$
 such that in these coordinates the
 system (\ref{bl}) is equivalent to the system 
\begin{equation}\label{bl1}
\frac{\partial^2 }{(\partial x^k)^2}\phi=F_k(x^k, \lambda_0,\lambda_1, ... ,\lambda_{n-1})\phi,  \  \ k=0,1,...,n-1, 
\end{equation}
where the function $F_k$ depends on the variable $x^k$ 
and on the parameters $\lambda_0,\lambda_1, ... ,\lambda_{n-1}$. 
Then $\phi$ is the  product 
$$X_0(x^0)X_1(x^1)...X_{n-1}(x^{n-1}),$$
 and each $X_k$ is a solution of the $k$-th equation of  (\ref{bl1}) so that
 we 
reduced 
 the system of partial differential equations (\ref{bl}) to the system of
ordinary differential equations 
\begin{equation} 
\nonumber
\frac{\partial^2 }{(\partial x^k)^2}X_k(x^k)=F_k(x^k, \lambda_0,\lambda_1, ... ,\lambda_{n-1})X_k(x^k),  \  \ k=0,1,...,n-1.
\end{equation}

\vspace{2ex}

\noindent {\bf The space of the metrics, geodesically equivalent to a given one: }

 \noindent  If two metrics are geodesically equivalent 
then there exist at least one-parametric
family of geodesically equivalent metrics and this family
includes these two metrics \cite{TM1} or \cite{BFM}.
It is possible to  show  that the set
of metrics, geodesically equivalent to a given metric, is a manifold. 
What is the dimension of this manifold?

  Even locally, there exist metrics that admit no (non-trivial) geodesically 
  equivalent \cite{Mikes1}.  

  Even  locally, the dimension of this space does not  exceed   
  $\frac{(n+1)(n+2)}{2}$ 
  and is equal to $\frac{(n+1)(n+2)}{2}$ 
  only for the metrics of constant curvature \cite{Mikes1}.

Let $g, \bar g$  be geodesically equivalent  on $M^n$.
Suppose $M^n$ is closed  and suppose the metrics are
 strictly non-proportional    almost everywhere.
Then the geodesic flow 
of the metric $g$ is completely  integrable, 
and almost all orbits 
lie at the corresponding Liouville tori. 
Suppose that  the geodesic flow
is  non-resonant.
 Then the Liouville foliation 
is uniquely  definite, and any integral of the geodesic flow commutes with the 
integrals $I_0,...,I_{n-1}$. 
Assume in addition, that there are sufficiently many caustics of the 
Liouville tori of the geodesic flow: almost each point of the surface is an 
intersection of $n$ caustics. Then the dimension of the space of the metrics (modulo multiplication by a constant), geodesically equivalent to 
the metric $g$, is equal one.

The geodesic flows of the metric of the
 ellipsoid and of the  Poisson sphere (for different $a_1, a_2, ..., a_n$)
satisfies all these conditions.

The authors are grateful to  Professors  A.V. Aminova,
W. Ballmann, V. Bangert,   A.V. Bolsinov, V. Cortes,  A. T. Fomenko, M. Gromov,
E. Heintze, A.S. Mistchenko,
  A.M. Perelomov,   S. Tabachnikov for interesting discussions.

\end{document}